\documentclass[12pt]{article}
\usepackage{amsmath,amsfonts,amssymb,amsthm}
\usepackage{bbm,mathrsfs,bm,mathtools}
\title{Pseudo-cones and measure transport}
\author{Rolf Schneider}
\date{}
\newcommand{\Sn}{{\mathbb S}^{n-1}}

\newcommand{\R}{{\mathbb R}}

\newcommand{\N}{{\mathbb N}}

\newcommand{\D}{{\rm d}}

\newtheorem{theorem}{Theorem}

\begin{document}
\maketitle

\begin{abstract}
{A recent result on the Gauss image problem for pseudo-cones can be interpreted as a measure transport, performed by the reverse radial Gauss map of a pseudo-cone. We find a cost function that is extremized by this transport map, and we prove an analogue of Rockafellar's characterization of the subdifferentials of convex functions.}\\[2mm]
\end{abstract}

\section{Introduction and formulation of results}\label{sec1}

A pseudo-cone $K\subset\R^n$ is a closed convex set not containing the origin $o$ and satisfying $\lambda K\subseteq K$ for $\lambda\ge 1$. In the following, we consider only pseudo-cones with a fixed recession cone. (Recall that the recession cone of $K$ is defined by ${\rm rec}\,K=\{z\in \R^n: K+z\subseteq K\}$; see Rockafellar \cite[Sect. 8]{Roc70}. For further details on convex sets, we refer to \cite{Sch14}.) We assume that a closed convex cone $C\subset \R^n$, pointed and with interior points, is given. A $C$-pseudo-cone $K$ is then a pseudo-cone with recession cone C.  A nonempty closed convex set $K\subset \R^n$ is a $C$-pseudo-cone if and only if $o\notin K$ and $K+C=K\subset C$.
Denoting by $\Sn$ the unit sphere of $\R^n$, we set
$$ \Omega_C:= \Sn \cap {\rm int}\,C,\qquad \Omega_{C^\circ}:= \Sn \cap{\rm int}\,C^\circ,$$
where $C^\circ:= \{x\in \R^n:\langle x,y\rangle \le 0\,\forall y\in C\}$ is the dual cone of $C$. We denote by $\langle \cdot\,,\cdot\rangle$ the scalar product, by int the interior, by cl the closure, and by bd the boundary. The set of all $C$-pseudo-cones in $\R^n$ is denoted by $ps(C)$.

Each set $ps(C)$ may be considered as a counterpart to the set of convex bodies containing the origin in the interior. There is a copolarity with properties similar to those of the polarity of convex bodies (see \cite{Sch24}), and there are Minkowski type problems for the (now possibly infinite) surface area measure and cone-volume measure and their generalizations and analogues.

Let $K\in ps(C)$. The radial function $\rho_K: \Omega_C\to (0,\infty)$ is defined by 
$$\rho_K(v):= \min\{r>0: rv\in K\}\quad\mbox{for } v\in \Omega_C.$$ 
By a normal vector of $K$ we always mean an outer unit normal vector. Each normal vector of $K$ belongs to ${\rm cl}\,\Omega_{C^\circ}$. Let $\omega_K$ be the set of all $u\in{\rm cl}\,\Omega_{C^\circ}$ that are normal vectors at more than one point of $K$. It is known that $\omega_K$ can be covered by countably many sets of finite $(n-2)$-dimensional Hausdorff measure and hence has Hausdorff dimension at most $n-2$. If $u\in \Omega_{C^\circ}\setminus \omega_K$, there is a unique vector $v\in\Omega_C$ such that $u$ is a normal vector of $K$ at $\rho_K(v)v$. (Note that if $u$ is attained at a point of $K\cap {\rm bd}\,C$, then $u\in\omega_K$.) We write $ v=\alpha_K^*(u)$ and call the map $\alpha_K^*:\Omega_{C^\circ}\setminus \omega_K\to \Omega_C$ thus defined the {\em reverse radial Gauss map} of $K$.

Starting point of this note is the following theorem, which was proved in \cite{Sch25}. Here we denote by $P(X)$ the set of Borel probability measures on a topological space $X$.

\vspace{2mm}

\noindent{\bf Theorem A.} {\em Let $\mu\in P(\Omega_{C^\circ})$ and $\nu\in P(\Omega_C)$ and suppose that $\mu$ is zero on sets of Hausdorff dimension $n-2$. Then there exists a $C$-pseudo-cone $K\in ps(C)$ such that} $(\alpha_K^*)\texttt{\#}\mu=\nu${\em, where $\alpha_K^*$ is the ($\mu$-almost everywhere on $\Omega_{C^\circ}$ defined) reverse radial Gauss map of $K$.}

\vspace{2mm}

Here $(\alpha_K^*)\texttt{\#}\mu=\nu$ means that $\alpha_K^*$ pushes $\mu$ forward to $\nu$, that is, $\mu((\alpha_K^*)^{-1}(\eta)) = \nu(\eta)$ for each Borel set $\eta\subseteq \Omega_C$.

We remark that $K$ in Theorem A is not uniquely determined; any dilate of $K$ has the same property. However, uniqueness up to dilatations has only been proved in \cite{Sch25} under additional assumptions.

Theorem A should be compared to the following  well-known result (we refer, e.g., to  McCann \cite{McC95}).

\vspace{2mm}

\noindent{\bf Theorem B.} (Brenier--McCann) {\em Let $\mu,\nu\in P(\R^n)$ and suppose that $\mu$ is zero on sets of Hausdorff dimension $n-1$. Then there exists a convex function $f:\R^n\to (-\infty,\infty]$ such that} $(\nabla f)\texttt{\#}\mu=\nu${\em, where $\nabla f$ denotes the ($\mu$-almost everywhere on ${\rm dom}\,f$ defined) gradient of $f$.}

\vspace{2mm}

Thus, to the gradient map of a convex function in Theorem B there corresponds in Theorem A the reverse radial Gauss map of a $C$-pseudo-cone. In either case, the theorem provides a map that transports the measure $\mu$ to the measure $\nu$. Usually in measure transportation theory (we refer to Lecture II of \cite{Art23} and also to the books \cite{RR99}, \cite{Vil03}, \cite{Vil09}), one is interested in transport plans or maps that extremize a certain total cost. Although Theorem A was proved by a different method, it implies that the obtained transportation map minimizes a certain total cost, namely that of the cost defined by
\begin{equation}\label{1.1}
c(u,v):= \log|\langle u,v\rangle|,\quad (u,v)\in \Omega_{C^\circ}\times\Omega_C.
\end{equation}
This is shown by the following theorem. Here we denote (for $\mu,\nu$ as in Theorem A) by ${\mathcal T}$ the set of all measurable, $\mu$-almost everywhere defined mappings $T$ from $\Omega_{C^\circ}$ to $\Omega_C$ with $T\texttt{\#}\mu=\nu$.

\begin{theorem}\label{T1.1}
If $\mu,\nu,K$ are as in Theorem A, then
$$ \int_{\Omega_{C^\circ}} c(u,\alpha_K^*(u))\,\mu(\D u) \le \int_{\Omega_{C^\circ}}c(u,T(u))\,\mu(\D u)$$
for all $T\in{\mathcal T}$.
\end{theorem}

Of course, this is trivial if the left side is equal to $-\infty$. But this is not always the case, for example, not if the support of $\mu$ is a closed subset of $\Omega_{C^\circ}$.

Theorem \ref{T1.1} was suggested by an investigation of Oliker \cite{Oli07}, who first treated Aleksandrov's integral curvature problem for convex bodies by a variational argument and then established a connection to optimal transport. There is an essential difference concerning the Gauss image problem for convex bodies and for pseudo-cones, when its relation to measure transport is considered (as done by Bertrand [2] for convex bodies). Taking the negatives of the cost functions in both cases, so that they become nonnegative, we see that the found transport map in the first case minimizes the total cost, whereas in the second case it maximizes it. For pseudo-cones $K$, this is not surprising, since $|\langle u, \alpha_K^*(u)\rangle|$ (which is less than $1$) is often close to zero, hence its negative logarithm is large. On the other hand, for special cones $C$, there are even transport maps (necessarily far away from a reverse radial Gauss map) for which the (nonnegative) total cost becomes zero. This is the case if $C$ is chosen as a circular cone of such size that the map defined by $T(u)=-u$ becomes a diffeomeorphism of $\Omega_{C^\circ}$ to $\Omega_C$, and the measures $\mu,\nu$ are the normalized restrictions of spherical Lebesgue measure on $\Omega_{C^\circ}$ and $\Omega_C$.

The gradients of the convex function $f$ appearing in Theorem B are subsumed in the subdifferential $\partial f$ of $f$, which is defined by
$$ \partial f:= \{(x,y)\in\R^n\times\R^n: f(x)<\infty\,\mbox{and}\, f(z)-f(x)\ge \langle y, z-x\rangle\,\forall z\in\R^n\}.$$
The subdifferentials of convex functions are characterized by Rockafellar's \cite{Roc66} classical theorem (see also \cite[Thm. 24.8]{Roc70} and \cite{PKW23}), which plays, together with its extensions, an essential role in measure transportation theory. For a general cost function $c:X\times Y\to\R$, where $X,Y$ are arbitrary sets, one says that a set $S\subset X\times Y$ is {\em $c$-cyclically monotone} if
$$ \sum_{i=1}^N c(x_i,y_i) \le \sum_{i=1}^N c(x_i,y_{\sigma(i)})$$
for all $n\in \N$, all $(x_i,y_i)\in S$ and all permutations $\sigma$ of $\{1,\dots,N\}$. For $X=Y=\R^n$, {\em cyclically monotone} means $c$-cyclically monotone for $c(x,y):=-\langle x,y\rangle$.

\vspace{2mm}

\noindent{\bf Theorem C.} (Rockafellar) {\em Let $S\subset \R^n\times\R^n$. There exists a convex function $f:\R^n\to (-\infty,\infty]$ with $S\subseteq \partial f$ if and only if $S$ is cyclically monotone.}

\vspace{2mm}

Since gradients of convex functions and reverse radial Gauss maps of $C$-pseudo-cones play analogous roles in Theorems B and A, the question arises whether there is a notion of subdifferential for pseudo-cones, leading to an analogue of Rockafellar's theorem. In fact, if we define the {\em pseudo-subdifferential} of $K\in ps(C)$ by
$$ \partial^\bullet K:= \{(v,u)\in\Omega_C\times\Omega_{C^\circ}: \mbox{$u$ is a normal vector of $K$ at $\rho_K(v)v$}\},$$
then 
$$ (\alpha_K^*(u),u)\in\partial^\bullet K\quad\mbox{for almost all $u\in\Omega_{C^\circ}$},$$
and the following theorem holds.

\begin{theorem}\label{T1.2}
Let $S\subset \Omega_C\times\Omega_{C^\circ}$. There exists a $C$-pseudo-cone $K\in ps(C)$ with $S\subseteq \partial^\bullet K$ if and only if $S$ is c-cyclically monotone for the cost function $c$ given by $c(v,u)=\log|\langle v,u\rangle|$.
\end{theorem}

We prove Theorem \ref{T1.1} in Section \ref{sec2} and Theorem \ref{T1.2} in Section \ref{sec3}.

\section{Proof of Theorem \ref{T1.1}}\label{sec2}

Let $\mu,\nu$ be as in Theorem A, let ${\mathcal T}$ be defined as before Theorem \ref{T1.1}, and let $K\in ps(C)$. The support function of $K$ is defined by
$$ h_K(u):= \sup\{\langle u,y\rangle:y\in K\}\quad\mbox{for } u\in{\rm cl}\,\Omega_{C^\circ}.$$
Since $h_K\le 0$, we write $\overline h_K=-h_K$. By the definition of the support function we have
\begin{equation}\label{2.1}
\overline h_K(u)\le |\langle u,\rho_K(v)v\rangle|\quad\mbox{for }(u,v)\in \Omega_{C^\circ}\times \Omega_C.
\end{equation}
Here equality holds if $u$ is a normal vector of $K$ at $\rho_K(v)v$, therefore
\begin{equation}\label{2.2}
\overline h_K(u) =|\langle u,\alpha_K^*(u)\rangle|\rho_K(\alpha_K^*(u))\quad\mbox{for } u\in\Omega_C^\circ\setminus\omega_K.
\end{equation}
From (\ref{2.1}) we get
$$ \log \overline h_K(u) -\log\rho_K(v) \le \log|\langle u,v\rangle|= c(u,v),$$
where $c$ is defined by (\ref{1.1}). For $T\in{\mathcal T}$ this gives
$$ \log\overline h_K(u)- \log \rho_K(T(u)) \le c(u,T(u))$$
for $\mu$-almost all $u\in \Omega_{C^\circ}$. Integration with the measure $\mu$ gives
$$ \int_{\Omega_{C^\circ}} \log\overline h_K(u) \,\mu(\D u) - \int_{\Omega_{C^\circ}} \log\rho_K(T(u))\,\mu(\D u) \le \int_{\Omega_{C^\circ}} c(u,T(u))\,\mu(\D u),$$
where at least the middle integral is greater than $-\infty$. Since $T\texttt{\#}\mu=\nu$ for $T\in{\mathcal T}$, the change of variables formula yields
\begin{equation}\label{2.3}
\int_{\Omega_{C^\circ}} \log\overline h_K(u) \,\mu(\D u) - \int_{\Omega_C} \log\rho_K(v)\,\nu(\D v) \le \int_{\Omega_{C^\circ}} c(u,T(u))\,\mu(\D u).
\end{equation}

Now let $K$ be a pseudo-cone as provided by Theorem A. Then $\alpha_K^*\in {\mathcal T}$, and the equality (\ref{2.2}) holds. Therefore, the inequality (\ref{2.3}) with $T=\alpha_K^*$ becomes an equality, that is,
$$ \int_{\Omega_{C^\circ}} \log\overline h_K(u) \,\mu(\D u) - \int_{\Omega_C} \log\rho_K(v)\,\nu(\D v) = \int_{\Omega_{C^\circ}} c(u,\alpha_K^*(u))\,\mu(\D u).$$
It follows that 
$$ \int_{\Omega_{C^\circ}} c(u,\alpha_K^*(u))\,\mu(\D u) \le \int_{\Omega_{C^\circ}} c(u,T(u))\,\mu(\D u)$$
for all $T\in{\mathcal T}$, as stated.

\section{Proof of Theorem \ref{T1.2}}\label{sec3}

Let $K\in ps(C)$. Let $m\in\N$ and $(v_i,u_i)\in\partial^\bullet K$ for $i=1,\dots,m$. It follows from the definition of the support function that
$$ h_K(u_i)=\langle v_i,u_i\rangle\rho_K(v_i)$$
and
$$ h_K(u_i) \ge \langle v_i,u_{\sigma(i)}\rangle\rho_K(v_i)$$
for all permutations $\sigma$ of $\{1,\dots,m\}$. Therefore
$$ \prod_{i=1}^m \langle v_i, u_i\rangle \prod_{i=1}^m \rho_K(v_i) = \prod_{i=1}^m h_K(u_i)  = \prod_{i=1}^m h_K(u_{\sigma(i)})\ge\prod_{i=1}^m \langle v_i,u_{\sigma(i)}\rangle \prod_{i=1}^m \rho_K(v_i)$$
and hence
$$ \prod_{i=1}^m \langle v_i,u_i\rangle \ge \prod_{i=1}^m \langle v_i,u_{\sigma(i)}\rangle,$$
equivalently (since $\langle v,u\rangle <0$ for $(v,u)\in\Omega_C\times\Omega_{C^\circ}$)
$$ \sum_{i=1}^m \log |\langle v_i,u_i\rangle| \le \sum_{i=1}^m \log |\langle v_i,u_{\sigma(i)}\rangle|$$
for all permutations $\sigma$ of $\{1,\dots,m\}$. We have shown that $\partial^\bullet K$ is $c$-cyclically monotone for the cost function $c$ defined by $c(v,u)=\log|\langle v,u\rangle|$.

Conversely, let $S\subset\Omega_C\times\Omega_{C^\circ}$ be any set that is $c$-cyclically monotone for $c(v,u)=\log|\langle v,u\rangle|$. To show that it satisfies $S\subseteq\partial^\bullet K$ for some $K\in ps(C)$, we use the generalization of Rockafellar's theorem to general cost functions, due to Rochet \cite{Roc87} and R\"uschendorf. Proofs (extending Rockafellar's argument) can be found in R\"uschendorf \cite[Lem. 2.1]{Rus96} and Rachev and R\"uschendorf \cite[Prop. 3.3.9]{RR99} (where we replace $c\to -c$ and $f\to-\varphi$); compare also \cite[Thm. 4.43]{Art23} and  \cite[Thm. 7]{PKW23}. We include into the following proposition parts of the proof, which we shall need.

We recall that in the general situation, where $X,Y$ are arbitrary sets and $c:X \times Y \to \R$ is a real cost function, the {\em $c$-subdifferential} of a function $\varphi: X \to [-\infty,\infty)$ is defined by
$$ \partial^c\varphi:= \{(x,y)\in X\times Y: c(x,y)-\varphi(x)\le c(z,y) -\varphi(z)\,\forall z\in X\}.$$

\vspace{2mm}

\noindent{\bf Proposition.}
{\em Let $X,Y$ be any sets and let $c:X\times Y\to\R$ be a cost function. If a set $S\subset X\times Y$ is $c$-cyclically monotone, then 
$S\subseteq\partial^{c}\varphi$, where $\varphi$ is the function defined by 
$$
\varphi(x) = \inf\Big\{c(x,y_m) - c(x_0,y_0) + \sum_{k=1}^m \left(c(x_k,y_{k-1}) -c(x_k,y_k)\right)\Big\}
$$
for $x\in X$, where $(x_0,y_0)\in S$ is arbitrarily chosen and where the infimum is over all $m\in\N$ and all $(x_i,y_i)\in S$, $i=1,\dots,m$.} We have $\varphi(x_0)=0$.

\vspace{2mm}

We apply this with $X=\Omega_C$, $Y=\Omega_{C^\circ}$, $ c(v,u)=\log|\langle v,u\rangle|$ and write $(v,u)$ for $(x,y)$. Let $S\subseteq \Omega_C\times\Omega_{C^\circ}$. Suppose that $S$ is $c$-cyclically monotone. With $(v_0,u_0)\in S$ arbitrarily chosen, we define $\varphi$ as above (note that $\varphi(v_0)=0$). The function $\varphi$ is an infimum of functions of the form
$$ f_i(v)= \log |\langle v, u_i\rangle|+a_i,\quad v\in\Omega_C,$$
with some $ u_i\in\Omega_{C^\circ}$ and some $a_i\in \R$, where $i$ is in some index set $I$. Writing $b_i:= e^{-a_i}$, we have
$$ f_i(v)= \log\frac{|\langle v,u_i\rangle|}{b_i} \quad\mbox{for }i\in I$$
and hence 
$$ \varphi(v) =\inf_{i\in I}  \log\frac{|\langle v,u_i\rangle|}{b_i}= \log \inf_{i\in I} \frac{|\langle v,u_i\rangle|}{b_i} =-\log\sup_{i\in I} \frac{b_i}{|\langle v, u_i\rangle|}.$$
With each function $f_i$, we associate the $C$-pseudo-cone
$$ K_i:= \{x\in C:\langle x,u_i\rangle \le -b_i\}.$$
From $\varphi(v_0)=0$ it follows that $f_i(v_0)\ge 0$ and hence $\langle v_0,u_i\rangle \le -b_i$, that is, $v_0\in K_i$. Since $\langle \rho_{K_i}(v)v,u_i\rangle =-b_i$ for $v\in\Omega_C$, the radial function of $K_i$ is given by
$$ \rho_{K_i}(v)= \frac{b_i}{|\langle v,u_i\rangle|},\quad v\in\Omega_C.$$
Define
$$ K:= \bigcap_{i\in I} K_i.$$
The intersection is not empty, since $v_0\in K$. Hence $K$ is a $C$-pseudo-cone, with radial function given by 
$$ \rho_K(v)= \sup_{i\in I} \frac{b_i}{|\langle v,u_i\rangle|},\quad v\in\Omega_C.$$
Now it follows  that 
$$ \varphi(v) = -\log\rho_K(v).$$
We have
\begin{eqnarray*}
(v,u)\in\partial^{c}\varphi &\Leftrightarrow& c(v,u)-\varphi(v) \le c(w,u) -\varphi(w)\quad\forall w\in\Omega_C\\
&\Leftrightarrow& \log|\langle v,u\rangle| +\log \rho_K(v) \le \log|\langle w,u\rangle| +\log \rho_K(w)\quad\forall w\in\Omega_C\\
&\Leftrightarrow& \langle v,u\rangle \rho_K(v) \ge \langle w,u\rangle \rho_K(w)\quad\forall w\in\Omega_C
\end{eqnarray*}
and
\begin{eqnarray*}
(v,u) \in \partial^\bullet K &\Leftrightarrow& \mbox{$u$ is a normal vector of $K$ at $\rho_K(v)v$}\\
&\Leftrightarrow& h_K(u) = \langle \rho_K(v)v,u\rangle\mbox{ and } h_K(u)\ge \langle z,u\rangle\quad\forall z\in K\\
&\Leftrightarrow& \langle v,u\rangle \rho_K(v) \ge \langle w,u\rangle \rho_K(w) \quad\forall w\in\Omega_C,
\end{eqnarray*}
since $h_K(u)=\max_{z\in\partial K}\langle z,u\rangle$ for $u\in\Omega_{C^\circ}$ and $h_K(u)\ge \langle z,u\rangle$ for all $z\in K$ if and only if $h_K(u)\ge \langle z,u\rangle$ for all $z\in{\rm bd}\,K$. Since $S\subseteq \partial^{c}\varphi$ by the Proposition, it follows that $S\subseteq\partial^\bullet K$. This completes the proof.

\noindent Author's address:\\[2mm]
Rolf Schneider\\Mathematisches Institut, Albert--Ludwigs-Universit{\"a}t\\D-79104 Freiburg i.~Br., Germany\\E-mail: rolf.schneider@math.uni-freiburg.de

\end{document}